\documentclass[11pt]{article}
\usepackage{amssymb,amsthm, url}
\usepackage{amsmath,color}
\usepackage{epsfig}

\newtheorem{theorem}{Theorem}[section]
\newtheorem{lemma}[theorem]{Lemma}

\newtheorem{example}[theorem]{Example}

\newcommand{\eoproof}{\hspace*{\fill} $\square$ \vspace{5pt}}

\newcommand{\R}{\mathbb R}
\newcommand{\Q}{\mathbb Q}
\newcommand{\Z}{\mathbb Z}

\DeclareMathOperator{\sat}{sat}
\DeclareMathOperator{\cone}{cone}
\DeclareMathOperator{\monoid}{semi-group}
\DeclareMathOperator{\lattice}{lattice}

\newcommand{\FourTiTwo}{{\tt 4ti2}}

\begin{document}

\title{Computing holes in semi-groups and its applications to transportation problems}

\author{Raymond Hemmecke, Akimichi Takemura, and Ruriko Yoshida}

\date{}

\maketitle

\begin{abstract}
An integer feasibility problem is a fundamental problem in many areas, 
such as operations research, number theory, and statistics.
To study a family of systems with no nonnegative integer solution, we focus 
on a commutative semigroup generated by a finite 
set of vectors in $\Z^d$ and its saturation. In this paper
we present an algorithm to compute an explicit description for the
set of holes which is the difference of a semi-group $Q$ generated
by the vectors and its saturation. We apply our procedure
to compute an infinite family of holes for the
semi-group of the $3\times 4\times 6$ transportation problem. Furthermore,
we give an upper bound for the entries of the holes when the set of
holes is finite. Finally, we present an algorithm to find all $Q$-minimal
saturation points of $Q$.
\end{abstract}

\section{Introduction}
The linear integer feasibility problem is to ask whether the system 
\begin{equation}\label{feas}
A x = b, \quad x \geq 0, 
\end{equation}
where $A \in \Z^{d\times n}$ and $b \in \Z^d$, has an integral solution or 
not.  In \cite{takemura-yoshida2006} we studied
a {\em generalized integer  feasibility problem}, that is,
to find %a finite representation of 
all  $b$ with no
nonnegative integral solution for a given $A$.
In recent years, the generalized integer linear feasibility problem has found 
applications 
in many research areas, such as number theory and statistics.
%\cite{takemura-yoshida2006}.
For example, in number theory, the {\em Frobenius problem} is to find the 
largest positive
integer $b$ such that there does not exist an integral solution in \eqref{feas}
with $d = 1$ (e.g.\ \cite{frob}, \cite{aardaletal3}).
In statistics, one can find an 
application in the data security problem of 
{\em multi-way contingency tables} \cite{dobra-karr-sanil2003}.
One of the challenge problems is the {\em 3-dimensional integer planar 
transportation problem} (3-DIPTP), that is, the question to decide whether the
set of \emph{integer} feasible solutions of the $r\times s\times
t$-transportation problem 
\[
\left\{x\in\Z^{rst}:\sum_{i=1}^r x_{ijk}=u_{jk},\sum_{j=1}^s x_{ijk}=v_{ik},\sum_{k=1}^t
x_{ijk}=w_{ij}, x_{ijk}\geq 0\right\}
\] 
is empty or not for a given right-hand sides $u,v,w$.
Vlach provides an excellent summary of attempts on 3-DIPTP \cite{Vlach:86}.
For sequential importance sampling \cite{chen-dinwoodie-sullivant}, 
non-existence of integral solution causes difficulties in its implementation.

Note that there exists a real nonnegative solution but there does not exist an integral 
nonnegative solution in  \eqref{feas} if and only if $b$ is in the
difference between the {\em semigroup} $Q$
generated by the column vectors of $A$ and its saturation 
% Therefore, we consider the difference of a semi-group $Q$
% generated by the columns of a matrix $A\in\Z^{d\times n}$ and its
% saturation 
$Q_{\sat}=\cone(A)\cap\lattice(A)$, where $\cone(A)$ is
the cone generated by the columns of $A$ and $\lattice(A)$ is the
lattice generated by the columns of $A$. We assume $\cone(A)$ to be
pointed.  
We call $H=Q_{\sat}\setminus Q$ the set of \emph{holes} of $Q$ and
call $Q$ \emph{normal} if $H=\emptyset$. $H$ may be finite or
infinite. 

In this paper, we present an algorithm which gives a finite description
of $H$. % even  when $H$ is infinite
Practically, even with all the currently available software
packages, checking normality of $Q$
%, that is whether $H$ is empty or not, 
is still a difficult computational question. Computing a finite
description of \emph{all} elements in $H$ is even more difficult.
The reader should note that for fixed matrix sizes $d$ and $n$,
there exists a \emph{polynomial size} encoding of the generating
function $f(H;z) = \sum_{h\in H} z^h$ (where
$z^h:=z_1^{h_1}\cdots z_d^{h_d}$) as a rational generating
function \cite{takemura-yoshida2006},\cite{bar}:
\[
f(H;z)=\sum_{i \in I} \gamma_i \frac{z^{\alpha_i}}{\prod_{j=1}^d
(1-z^{\beta_{ij}})}.
\]
Herein, $I$ is a finite (polynomial size) index set and all the
appearing data $\gamma_i\in\Q$ and $\alpha_i,\beta_{ij}\in\Z^d$ are
of size polynomial in the input size of $A$. In fact, this
observation is based on a result by Barvinok and Woods
\cite{Barvinok+Woods:2003}, who showed that there are such
\emph{short} rational function encodings for $Q$ and for $Q_{\sat}$,
and consequently, also for $f(H;z)=f(Q_{\sat};z)-f(Q;z)$. 
Although the proof by Barvinok and Woods is constructive, its
practical usefulness still has to be proven by an efficient
implementation. 
% Finally, note that once $f(H;z)$ has been computed,
% one can also decide in polynomial time (when the dimensions $d$ and
% $n$ are kept fix) whether the sum of rational functions encodes a
% polynomial or an infinite series. Thus, one can decide finiteness of
% $H$ \cite{takemura-yoshida2006}.
%
In contrast to the \emph{implicit}
representation via rational generating functions, in this paper, we present an
algorithm to compute an \emph{explicit} representation
of $H$, even for an infinite set $H$.  Such an explicit
representation need not be of polynomial size in the input size of
$A$. 
% Moreover, such an explicit representation acannot be recovered
% easily from an implicit short rational function encoding of
% $f(H;z)$.

This paper is organized as follows:  In Section \ref{mainthm} we 
set up basic notation and present our main results.  Section 
\ref{fundamental} shows a combinatorial algorithm to compute the set of
all {\em fundamental holes} of $Q$.  In Section \ref{hole} we describe
an algorithm to compute a finite representation of holes in $Q$.
Section \ref{trans} shows an application of the algorithm to 3-DIPTP and
in Section \ref{boundary} we describe the bounds on the size of entries in 
each hole in $Q$.  Finally in Section \ref{minimal} we show an algorithm to 
compute the set of all
{\sl $Q$-minimal saturation points}.

\section{Basic notation and the main results}\label{mainthm}

The main result in this paper is the following.
\begin{theorem}
\label{therem:main}
There exists an algorithm that computes for an integral matrix $A$ a
finite explicit representation for the set $H$ of holes of the
semigroup $Q$ generated by the columns of $A$, that is, the algorithm
computes (finitely many) vectors $h_i\in\Z^n$ and monoids $M_i$, each given
by a finite set of generators in $\Z^n$, $i\in I$, such that
\[
H=\bigcup_{i\in I}\; \left(\{h_i\}+M_i\right).
\]
\end{theorem}

In fact, we explicitly present such an algorithm. Note that $M_i$
could be trivial, that is, $M_i=\{0\}$. See Example \ref{Example 1}
below for an example of such an explicit representation.

One basic object needed in our construction is the set $F$ of
%$Q$-minimal holes $h\in H$, which we call 
fundamental holes.
We call $h\in H$ {\em fundamental} if there is no other hole $h'\in H$
such that $h-h'\in Q$. Note that in contrast to $H$, $F$ is always
finite. For every hole $h\in H$ there exists a
fundamental hole $f\in F$ such that $h\in f+Q$.
In view of these facts our algorithm consists of the following two main steps:
\begin{enumerate}
\item First, compute the set $F$ of fundamental holes.
\item Then, for each of the finitely many $f\in F$, compute an
explicit representation of the holes in $f+Q$.
\end{enumerate}

% We perform the second step of this algorithm for {\color{blue}
% the set of holes for the semi-group $Q$ of the $3\times 4\times 6$
% transportation problem in \cite{Vlach:86}}. We show that the set of 
% holes of $Q$
% belonging to $f+Q$ is \emph{not finite} by computing an explicit
% finite description for all these holes.

Moreover, we can use our algorithm to bound the entries of $h\in H$
in case that $H$ is finite.

{\bf Theorem  \ref{Lemma: Bounds}.} {\it Let $A\subseteq\Z^{d\times n}$
be of full row-rank. %and let $H$ be defined as above. 
Let $D(A)$
denote the maximum absolute value of the determinants of a $d\times
d$ submatrix of $A$. Moreover, let
$M_F(A)=\max_{i=1,\ldots,d}\sum_{j=1}^n |A_{ij}|-1$. Then, if $H$ is
finite, the inequality
\[
\|h\|_{\infty}\leq (d+1)M_F^2(A)D(A)
\]
holds for every $h\in H$. }

Finally, we can use the above approach to compute the set of all
$Q$-minimal saturation points of $Q$ (Section \ref{minimal}). 
Herein, we call $s\in Q$ a
\emph{saturation point} of $Q$, if $s+Q_{\sat}\subseteq Q$. The set
of all saturation points of $Q$ is denoted by $S$. We call $s\in S$
a $Q$-minimal saturation point if there is no other $s'\in S$ with
$s-s'\in Q$.  The set $S$ of saturation points, considered as a $Q$-module,
is often called a {\em conductor ideal} (e.g.\  \cite{BGT}). $S$ is
finitely generated as a $Q$-module and hence the set of $Q$-minimal
saturation points is finite.

We illustrate the above notions with the following simple example.

\begin{example}\label{Example 1}\rm
Let
\[
A=\left(\begin{array}{cccc} 1 & 1 & 1 & 1\\
0 & 2 & 3 & 4\\
\end{array}\right)
\]
with the single fundamental hole $(1,1)^\intercal$ and with
infinitely many holes
\[
H=\{(1,1)^\intercal+\alpha\cdot (1,0)^\intercal:\alpha\in\Z_+\},
\]
where $\Z_+$ denote the set of nonnegative integers. 
$Q$ has three $Q$-minimal
saturation points: $(1,2)^\intercal$, $(1,3)^\intercal$, and
$(1,4)^\intercal$. %\eoproof
See Figure
\ref{Figure: Example 1}. \eoproof

\begin{figure}[tbh]
\begin{center}
\epsfig{file=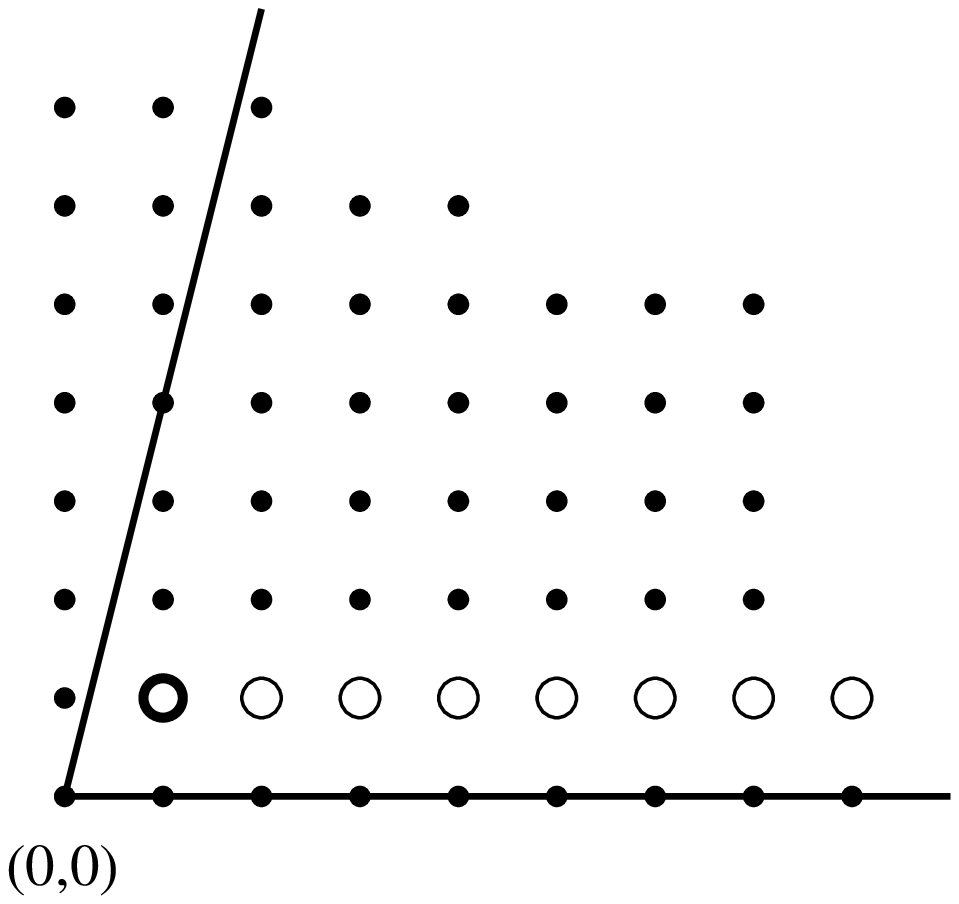, height=5cm}\\[.3cm]
\refstepcounter{figure} \label{Figure: Example 1} Figure \thefigure:
Non-holes, holes and fundamental hole for Example \ref{Example 1}
\end{center}
\vspace*{-0.3cm}
\end{figure}
\end{example}

% {\bf Example \ref{Example 1} cont.} 

In the following two sections we demonstrate how to perform the steps
for Theorem  \ref{therem:main}
algorithmically. We accompany the theoretical construction with our
running example, Example \ref{Example 1}.

\section{Computing the fundamental holes $F$}\label{fundamental}

In this section, we compute the set $F$ of fundamental 
%(that is                              %$Q$-minimal) 
holes $h\in H$. To enumerate all fundamental holes, we
first give a short proof that the number of fundamental holes is
indeed finite. Let $A_{.i}$ denote the $i$th column of $A$.
\begin{lemma}[Takemura and Yoshida \cite{takemura-yoshida2006}]\label{holes}
The set $F$ of fundamental holes is a subset of
\[
P:=\left\{\sum_{i=1}^n \lambda_i
A_{.i}:0\leq\lambda_1,\ldots,\lambda_n<1\right\} \cap \Z^d.
\]
\end{lemma}
\begin{proof}
Each $f\in F$ lies in $\cone(A)$ and thus can be written as
$f=A\lambda=\sum_{i=1}^n \lambda_i A_{.i}$ for some $\lambda\geq 0$.
If $\lambda_j\geq 1$ for some $j\in\{1,\ldots,n\}$ then
$f'=f-A_{.j}\in\cone(A)\cap\lattice(A)$ would contradict the
$Q$-minimality of $f$, since $f-f'=A_{.j}\in Q$. Consequently,
$\lambda_i<1$ for all $i$.
\end{proof}
This shows that $F$ is finite and also gives a finite procedure to
enumerate $F$:
\begin{itemize}
\item Enumerate $P\cap\lattice(A)$.
\item Check for each $z\in P\cap\lattice(A)$ whether $z$ is a
fundamental hole or not by solving $A\lambda=z, \lambda\in\Z^n_+$
and by checking whether $z-A_{.i}\in P\cap\lattice(A)\subseteq Q$
for some $i$.
\end{itemize}

Practically, this construction can be sped-up as follows. First
compute the (unique) minimal Hilbert basis (or better: integral
basis) $B$ of $\cone(A)\cap\lattice(A)$. Again, similarly as above,
one can show that $B\subseteq P$. If $B$ contains no hole of $Q$,
$Q$ must be normal. Otherwise, every hole of $Q$ appearing in $B$
must be fundamental, since $B$ is minimal. Finally, if $f\in F$ is
not in $B$, $f$ can be written as a nonnegative integer linear
combination of elements in $B$, since $f\in\cone(A)\cap\lattice(A)$
and since $B$ is a Hilbert basis (integral basis) of
$\cone(A)\cap\lattice(A)$. This representation cannot have summands
that are not fundamental holes, since otherwise $f$ is would not be
fundamental. To see this, let
\[
f=\sum_{b\in B\cap F} \lambda_b b+ \sum_{b\not\in B\cap F} \mu_b
b,\;\;\; \lambda_b,\mu_b\in\Z_+\;\forall b,
\]
with $\sum_{b\not\in B\cap F} \mu_b b\neq 0$.
Observe, that
\[
f'=\sum_{b\in B\cap F} \lambda_b b
\]
must be a hole of $Q$, as otherwise $f$ is not a hole. But since
\[
f-f'=\sum_{b\not\in B\cap F} \mu_b b\in Q,
\]
$f$ cannot be a fundamental hole.

Thus  we can enumerate $F$ as follows:
\begin{itemize}
\item Compute the Hilbert basis (integral basis) $B$ of
$\cone(A)\cap\lattice(A)$.
\item Check whether each $z\in B$ is a fundamental hole or not,
that is, compute $B\cap F$.
\item Generate all nonnegative integer combinations of elements in $B\cap
F$ that lie in $P$ and check for each such $z$ whether it is a
fundamental hole or not.
\end{itemize}

{\bf Example \ref{Example 1} cont.} In our example, the lattice $L$
generate by the columns of $A$ is simply $\lattice(A)=\Z^2$. With
this, the Hilbert basis $B$ of $\cone(A)\cap\lattice(A)$ consists of
$5$ elements:
\[
B=\{(1,0)^\intercal,(1,1)^\intercal,(1,2)^\intercal,(1,3)^\intercal,(1,4)^\intercal\},
\]
out of which only $(1,1)^\intercal$ is a hole. Being in $B$,
$(1,1)^\intercal$ must be a fundamental hole. Thus, $B\cap
F=\{(1,1)^\intercal\}$. Constructing nonnegative integer linear
combinations of elements from $B\cap F$, we already see that the
combination $2\cdot (1,1)^\intercal=(2,2)^\intercal$ is an element
of $Q$ and consequently, there is no other fundamental hole in $Q$,
i.e.\ $F=\{(1,1)^\intercal\}$. \eoproof

\section{Computing the holes in $f+Q$}\label{hole}

In this section we discuss how to compute the holes in $f+Q$ for
each fundamental hole $f\in F$. Note that a point $z\in f+Q$ is
either a hole or it belongs to $Q$. That is, every non-hole in $f+Q$
belongs to $(f+Q)\cap Q$. Moreover, if $z\in (f+Q)\cap Q$ then also
$z+A\lambda\in (f+Q)\cap Q$ for all $\lambda\in\Z_+^n$. Thus we define
a monomial ideal $I_{A,f}\in\Q[x_1,\ldots,x_n]$ by
% generated by the monomials
\begin{equation}
\label{eq:Iaf}
I_{A,f}=\langle x^\lambda:\lambda\in\Z^n_+,f+A\lambda\in (f+Q)\cap
Q\rangle.
\end{equation}
By construction, $f+A\lambda$, $\lambda\in\Z^n_+$, is not a hole of
$Q$ if and only if $x^\lambda\in I_{A,f}$. Therefore, we are looking
for an explicit description of the monomials \emph{not} belonging to
the monomial ideal $I_{A,f}$. These monomials are usually called
\emph{standard monomials} and there are algorithms to compute an
explicit disjoint or non-disjoint representation of them once ideal
generators for $I_{A,f}$ are known. Via the (typically
non-injective) linear transformation $x^\lambda\mapsto f+A\lambda$,
one recovers an explicit (usually non-disjoint) representation of
all holes of $Q$ in $f+Q$.

It remains to find (minimal) generators for $I_{A,f}$. The minimal
generators correspond to the $\leq$-minimal elements in the set
\[
L_{A,f}=\{\lambda\in\Z^n_+:\exists\mu\in\Z^n_+\text{ such that }
f+A\lambda=A\mu\}.
\]
To compute these minimal elements directly inside this projection is
a hard computational task and deserves further investigation. Let us
therefore compute a usually \emph{non-minimal} generating set for
$I_{A,f}$ from a higher-dimensional problem.

\begin{lemma}
Let $M$ be the set of $\leq$-minimal solutions $(\lambda,\mu)$ to
$f+A\lambda=A\mu$, $(\lambda,\mu)\in\Z^{2n}_+$. Then
\[
I_{A,f}=\langle x^\lambda:\exists\mu\in\Z^n_+\text{ such that }
(\lambda,\mu)\in M\rangle.
\]
\end{lemma}

\begin{proof}
%\boproof 
Let $\lambda_0\in L_{A,f}$ be $\leq$-minimal. We show now
that there exists some $\mu_0\in\Z^n_+$ such that
$(\lambda_0,\mu_0)$ is a $\leq$-minimal solution to
$f+A\lambda=A\mu$, $(\lambda,\mu)\in\Z^{2n}_+$. Then, as claimed,
the minimal generator $x^{\lambda_0}$ is contained in the given set
of generators for $I_{A,f}$.

Suppose on the contrary, that for every $\mu\in\Z^n_+$ the vector
$(\lambda_0,\mu)$ is \emph{not} a $\leq$-minimal solution to
$f+A\lambda=A\mu$, $(\lambda,\mu)\in\Z^{2n}_+$. Let $\mu_0$ be a
$\leq$-minimal solution to $f+A\lambda_0=A\mu$, $\mu\in\Z^{n}_+$.
Then, by our assumption, there is some vector
$(\lambda',\mu')\in\Z^{2n}_+$ with $f+A\lambda'=A\mu'$,
$(\lambda',\mu')\leq (\lambda_0,\mu_0)$, and $(\lambda',\mu')\neq
(\lambda_0,\mu_0)$. If $\lambda'\neq\lambda_0$ holds, we have a
contradiction to $\lambda_0$ being $\leq$-minimal in $L_{A,f}$. If
$\lambda'=\lambda_0$ and $\mu'\neq\mu_0$ holds, we have a
contradiction to $\mu_0$ being a $\leq$-minimal solution to
$f+A\lambda_0=A\mu$, $\mu\in\Z^{n}_+$. This shows that
$(\lambda_0,\mu_0)$ is a $\leq$-minimal solution to
$f+A\lambda=A\mu$, $(\lambda,\mu)\in\Z^{2n}_+$, as we wanted to
show. 
\end{proof}
%\eoproof

{\bf Example \ref{Example 1} cont.} Let $f=(1,1)^\intercal$ and
consider $(f+Q)\cap Q$. The linear system to solve is
\[
\begin{array}{rcrcrcrcrcrcrcrcr}
1 & + & \lambda_1 & + & \lambda_2 & + & \lambda_3 & + & \lambda_4
& = & \mu_1 & + & \mu_2 & + & \mu_3 & + & \mu_4\\
1 & & & + & 2\lambda_2 & + & 3\lambda_3 & + & 4\lambda_4
& = & &   & 2\mu_2 & + & 3\mu_3 & + & 4\mu_4\\
\end{array}
\]
with $\lambda_i,\mu_j\in\Z_+, i,j\in\{1,2,3,4\}$.

\FourTiTwo{} gives the following $5$ minimal inhomogeneous solutions
$(\lambda,\mu)\in\Z^8_+$:
{\footnotesize
\[
(0,0,0,2,0,0,3,0)^\intercal, (0,1,0,0,1,0,1,0)^\intercal,
(0,0,1,0,1,0,0,1)^\intercal (0,0,1,0,0,2,0,0)^\intercal,
(0,0,0,1,0,1,1,0)^\intercal.
\]
}
Thus, we get the monomial ideal
\[
I_{A,f}=\langle x_4^2,x_2,x_3,x_3,x_4\rangle=\langle
x_2,x_3,x_4\rangle,
\]
whose set of standard monomials is $\{x_1^\alpha:\alpha\in\Z_+\}$.
Thus, the set of holes in $f+Q$ is explicitly given by
\[
\{f+\alpha A_{.1}:\alpha\in\Z_+\}=\{(1,1)^\intercal+\alpha
(1,0)^\intercal:\alpha\in\Z_+\},
\]
as already claimed above. \eoproof

\section{Infinitely many holes for $3\times 4\times 6$
transportation problem}\label{trans}

In this section, we apply the procedure from the last section to the
semi-group $Q$ spanned by the matrix $A$ defining a $3\times 4\times
6$ transportation problem. Already in 1986, Vlach \cite{Vlach:86}
has shown that this semi-group is not normal by explicitly stating a
hole $f$, which is fundamental. He showed that $Ax=f, x\geq 0$ has a
(unique) rational solution, which in turn is fractional. In the
following, we construct a finite representation of all holes of $Q$
belonging to $f+Q$. We show that there are in fact \emph{infinitely
many} such holes.

The semi-group of the $3\times 4\times 6$ transportation problem is
of special interest, as it is the smallest
three-dimensional transportation problem for which it is known that the
associated semi-group is not normal. The only two cases of (also
higher-dimensional) transportation problems for which the normality
question is still open are those of sizes $3\times 4\times 5$ and
$3\times 5\times 5$ \cite{Ohsugi+Hibi:2007}. The previously open case
$3\times 4\times 4$ 
has been solved by the third author using the software package
NORMALIZ \cite{brunskoch}. The associated semi-group is normal.

For the $3\times 4\times 6$ problem, 
a fundamental hole $f$ given by Vlach \cite{Vlach:86}
% In {\color{blue}the set of fundamental holes}, the right-hand side $f$ for the counts
% along the coordinate axes 
is defined by the following three matrices:
\[
\left(
\begin{array}{ccc}
1 & 1 & 1\\
1 & 1 & 1\\
1 & 1 & 1\\
1 & 1 & 1\\
\end{array}
\right),
\left(
\begin{array}{ccc}
1 & 1 & 0\\
1 & 1 & 0\\
1 & 0 & 1\\
1 & 0 & 1\\
0 & 1 & 1\\
0 & 1 & 1\\
\end{array}
\right), \text{ and }
\left(
\begin{array}{cccc}
1 & 0 & 0 & 1\\
0 & 1 & 1 & 0\\
1 & 1 & 0 & 0\\
0 & 0 & 1 & 1\\
1 & 0 & 1 & 0\\
0 & 1 & 0 & 1\\
\end{array}
\right).
\]
The unique point in this $3\times 4\times 6$ transportation polytope
$\{z\in\R^{72}:Az=f, z\geq 0\}$ is
\[
z^*=\frac{1}{2}\left(
\begin{array}{c|c|c|c|c|c}
\begin{array}{ccc}
1 & 1 & 0 \\
0 & 0 & 0 \\
0 & 0 & 0 \\
1 & 1 & 0 \\
\end{array} &
\begin{array}{ccc}
0 & 0 & 0 \\
1 & 1 & 0 \\
1 & 1 & 0 \\
0 & 0 & 0 \\
\end{array} &
\begin{array}{ccc}
1 & 0 & 1 \\
1 & 0 & 1 \\
0 & 0 & 0 \\
0 & 0 & 0 \\
\end{array} &
\begin{array}{ccc}
0 & 0 & 0 \\
0 & 0 & 0 \\
1 & 0 & 1 \\
1 & 0 & 1 \\
\end{array} &
\begin{array}{ccc}
0 & 1 & 1 \\
0 & 0 & 0 \\
0 & 1 & 1 \\
0 & 0 & 0 \\
\end{array} &
\begin{array}{ccc}
0 & 0 & 0 \\
0 & 1 & 1 \\
0 & 0 & 0 \\
0 & 1 & 1\\
\end{array}
\end{array}
\right).
\]
If $A_{.,ijk}$ denotes the column of $A$ corresponding to variable
$z_{ijk}$, then the monomial ideal $I_{A,f}$ constructed in the
previous section is generated by the $48$ monomials $x_{ijk}$ for
which $z^*_{ijk}=0$. This can be shown as follows.

Firstly, $1\not\in I_{A,f}$, since $f\not\in Q$. Secondly, using {\tt
4ti2}, one verifies for each of these $48$ indices that
$f+A_{.,ijk}=A\mu$ has a nonnegative integer solution
$\mu\in\Z^{72}$, by explicitly constructing such a solution. Finally,
as it remains to look only for ideal generators of $I_{A,f}$ not
divisible by the $48$ monomials $x_{ijk}$ for which $z^*_{ijk}=0$, the
linear system from the previous section simplifies to
\[
f+A'\lambda'=A\mu,\lambda\in\Z^{24}_+,\mu\in\Z^{72}_+,
\]
where $A'$ is formed out of the $24$ columns $A_{.,ijk}$ of $A$ for
which $z^*_{ijk}>0$. This system does not have an integral solution.
In fact, the only real solution is $(\lambda',\mu)=(0,z^*)$. To see
this one either solves this system, for example using {\tt 4ti2}, or
one observes that the vector $f+A'\lambda'$ has many zero entries
that are present for arbitrary choices of $\lambda'$. These zero
entries imply that $\mu_{ijk}=0$ for all triples $ijk$ for which
$z^*_{ijk}=0$. The remaining linear system
\[
f+A'\lambda'=A\mu',\lambda\in\Z^{24}_+,\mu'\in\Z^{24}_+,
\]
has a unique real solution, namely $(0,{z^*}')$, which can be
checked by applying Gaussian elimination to the system
$A(\mu'-\lambda')=f, (\mu'-\lambda')\in\R^{24}$.

Thus, the set of holes of $Q$ belonging to $f+Q$ are given by
$f+\monoid(A')$.

\section{Computing bounds}\label{boundary}

For this section, let us assume that the set $H$ is finite. We will
now use our approach above to establish bounds on the size of the
entries for each $h\in H$. Clearly, such a bound can then be used to
show that $H$ cannot be finite if a hole with sufficiently big
entries has been found.

\begin{theorem}\label{Lemma: Bounds}
Let $A\subseteq\Z^{d\times n}$ be of full row-rank. 
% and let $H$ be defined as above. 
Let $D(A)$ denote the maximum absolute value of
the determinants of a $d\times d$ submatrix of $A$. Moreover, let
$M_F(A)=\max_{i=1,\ldots,d}\sum_{j=1}^n |A_{ij}|-1$. Then, if $H$ is
finite, the inequality
\[
\|h\|_{\infty}\leq (d+1)M_F^2(A)D(A)
\]
holds for every $h\in H$.
\end{theorem}

\begin{proof} First, we can bound the elements $f\in F$ using the
relation
\[
F\subseteq \left\{\sum_{j=1}^n \lambda_j
A_{.j}:0\leq\lambda_1,\ldots,\lambda_n<1\right\}.
\]
Thus,
\[
|f^{(i)}|\leq \sum_{j=1}^n
|A_{ij}|-1\leq\max_{i=1,\ldots,d}\sum_{j=1}^n |A_{ij}|-1=:M_F(A)
\]
holds for all $f\in F$ and all $i=1,\ldots,d$.

Next, as $H$ is finite, all ideals $I_{A,f}$, $f\in F$, must have a
finite set of standard pairs, which is equivalent to saying that
there must be a monomial generator $x_i^{\alpha_i}$ for every
$i=1,\ldots,n$. Such a monomial generator corresponds to a minimal
inhomogeneous solution $(\alpha_i,\mu)$ to $f+\alpha_i A_{.i}=A\mu,
\alpha_i\in\Z_+,\mu\in\Z^n_+$. Let us now bound the values for such
a minimal $\alpha_i$.

First, the minimal inhomogeneous solutions $(\alpha_i,\mu)$ to
$f+\alpha_i A_{.i}=A\mu, \alpha_i\in\Z_+,\mu\in\Z^n_+$ correspond
exactly to the minimal homogeneous solutions to $fu+\alpha_i
A_{.i}-A\mu=0, \alpha_i,u\in\Z_+,\mu\in\Z^n_+$ with $u=1$. Each
entry in a minimal homogeneous solutions, however, can be bounded by
$(d+1)$ times the maximum absolute value $D(f\ A_{.i}\-A)$ of the
determinants of a maximal submatrix of the defining matrix $(f\
A_{.i}\-A)$.

Thus, in particular,
\[
\alpha_i\leq (d+1) D(f\ A_{.i}\-A)\leq (d+1)\max_{j=1,\ldots,d}
|f^{(j)}|\cdot D(A_{.i}\ -A)=(d+1)M_F(A)\cdot D(A).
\]
Consequently, we can bound the entries of a hole in $f+Q$ by giving
bounds for
\[
f+\sum_{j=1}^n (\alpha_j-1)A_{.j}.
\]
For $h\in (f+Q)\cap H$, the $i$th entry is bounded as
\begin{eqnarray*}
h^{(i)}& \leq & |f^{(i)}|+\sum_{j=1}^n (\alpha_j-1)|A_{ij}|\\
& \leq & M_F(A) + \sum_{j=1}^n ((d+1)M_F(A)D(A)-1)|A_{ij}|\\
& =    & M_F(A) + ((d+1)M_F(A)D(A)-1)\sum_{j=1}^n |A_{ij}|\\
& \leq & M_F(A) + ((d+1)M_F(A)D(A)-1)M_F(A)\\[15pt]
& =    & (d+1)M_F^2(A)D(A).
\end{eqnarray*}
As this bound is independent on $f\in F$, we have
\[
\|h\|_{\infty}\leq (d+1)M_F^2(A)D(A)\;\;\;\forall h\in H,
\]
if $H$ is finite. 
\end{proof}
% \eoproof

{\bf Example \ref{Example 1} cont.} In our example, we have
\begin{itemize}
\item $d+1=3$,
\item $M_F(A)=\max(1+1+1+1,0+2+3+4)=9$, and
\item $D(A)=\max |2\times 2\text{ determinant of }A|=
|\det\left(\begin{smallmatrix}1 & 1\\0 & 4\\\end{smallmatrix}\right)|=4$.
\end{itemize}
Thus, if $H$ was finite, we would get the bound $\|h\|_{\infty}\leq
3\cdot 9^2\cdot 4=972$. In our example, however, one can easily
verify that $(1000,1)$ is a hole. Moreover, it violates the computed
bound. Consequently, $H$ cannot be finite. \eoproof

\section{Computing all $Q$-minimal saturation points}\label{minimal}

In this section, let $S$ denote the set of saturation points of $Q$,
that is, the set of all those $s\in Q$ such that
$s+Q_{\sat}\subseteq Q$. Let us now show how the above approach can
be used in order to compute $\min(S;Q)$, the set of all $Q$-minimal
points in $S$.  We also recover the known fact that $\min(S;Q)$ is
always finite.  We state the following theorem.

\begin{theorem}
\label{theorem:7.1}
\[
S=\bigcap_{f\in F} [((f+Q)\cap Q)-f]
\]
and  hence
\[
S = \{ A\lambda \mid x^\lambda\in\bigcap_{f\in F} I_{A,f} \},
\]
where $I_{A,f}$ is defined in (\ref{eq:Iaf}).
\end{theorem}

\begin{proof}
\begin{eqnarray*}
s\in S & \Leftrightarrow & s\in Q\text{ and }s+Q_{\sat}\subseteq
Q\;\;\;\text{(by
definition)}\\
& \Leftrightarrow & s\in Q\text{ and }s+H\subseteq Q\;\;\;
\text{(since }Q_{\sat}=Q\cup H\text{ and }s+Q\subseteq Q,\;\forall s\in Q\text{)}\\
& \Leftrightarrow & s\in Q\text{ and }s+F\subseteq Q\;\;\;
\text{(since }H\subseteq F+Q\text{)}\\
& \Leftrightarrow & s+f\in f+Q\text{ and }s+f\subseteq
Q\;\;\;\forall f\in F\\
& \Leftrightarrow & s+f\in (f+Q)\cap Q\;\;\;\forall f\in F.
\end{eqnarray*}
Consequently, we have
$
s\in S\Leftrightarrow s\in\bigcap_{f\in F} [((f+Q)\cap Q)-f]
$.
Furthermore with $s=A\lambda$ for some $\lambda\in\Z^n_+$ (as $s\in
Q$), we get
$
s\in S\Leftrightarrow x^\lambda\in\bigcap_{f\in F} I_{A,f}.
$
\end{proof}

Define 
\[
I_A=\bigcap_{f\in F} I_{A,f} .
\]
Then  $I_A$ is a monomial ideal being the intersection of the
monomial ideals $I_{A,f}$. $I_A$ can be found algorithmically, for
example with the help of Gr\"obner bases. The elements
$s\in\min(S;Q)$ correspond exactly to the minimal ideal generators
$x^\lambda$ of $I_A$ via the relation $s=A\lambda$. (Note, however,
that this relation need not be one-to-one. There may be many minimal
ideal generators corresponding to the same $Q$-minimal saturation
point.)

{\bf Example \ref{Example 1} cont.} In our example, we have
$I_A=I_{A,f}=\langle x_2,x_3,x_4\rangle$, as there exists only one
fundamental hole $f$. The three generators of $I_A$ correspond to
the three $Q$-minimal saturation points $(1,2)^\intercal$,
$(1,3)^\intercal$, and $(1,4)^\intercal$. \eoproof


\begin{thebibliography}{10}
\bibitem{4ti2}
4ti2 team. {\tt 4ti2}--Software package for algebraic, geometric and
combinatorial problems on linear spaces. Available at \url{http://www.4ti2.de/}.

\bibitem{aardaletal3}
K.~Aardal and A.~K.~Lenstra. Hard equality constrained integer knapsacks. In
{\sl Integer programming and combinatorial optimization}. 
Lecture Notes in Comput. Sci. Springer. (2002), 350--366.

\bibitem{frob}
J.~R.~Alfonsin, {\it The Diophantine Frobenius Problem}.
Oxford Lecture Series in Mathematics and Its Applications,
Oxford University Press, New York, 2006.


\bibitem{bar}
A.~I.~Barvinok. Polynomial time algorithm for counting integral
points in polyhedra when the dimension is fixed. {\sl Math. of
Operations Research} \textbf{19} (1994), 769--779.

\bibitem{Barvinok+Woods:2003}
A.~I.~Barvinok and K.~Woods. Short rational generating functions for
lattice point problems. {\sl J. Amer. Math. Soc.} \textbf{16}
(2003), 957--979.

\bibitem{BGT}
W.~Bruns, J.~Gubeladze and N.~V.~Trung.
Problems and algorithms for affine semigroups. {\it Semigroup Forum}
{\bf64} (2002), 180--212.

\bibitem{brunskoch}
W.~Bruns and R.~Koch. {\tt NORMALIZ},  computing normalizations of
affine semigroups, Available from \url{ftp://ftp.mathematik.uni-onabrueck.de/ pub/osm/kommalg/software/}.


\bibitem{chen-dinwoodie-sullivant}
Y.~Chen, I.~Dinwoodie,  and S.~Sullivant,
{\it Annals of Statistics}, \textbf{34} 
(2006), 523--545.

\bibitem{dobra-karr-sanil2003}
A.~Dobra, A.~F.~Karr, and P.~A.~Sanil. Preserving confidentiality of high-dimensional tabulated data: statistical and computational issues. {\sl Stat. Comput.}
\textbf{13} (2003), 363--370.

\bibitem{Ohsugi+Hibi:2007}
H.~Ohsugi and T.~Hibi. Toric ideals arising from contingency
tables. In {\sl Ramanujan Mathematical Society
Lecture Note Series}. No.4, In press. (2008).

\bibitem{takemura-yoshida2006}
A.~Takemura and R.~Yoshida. A generalization of the integer linear
infeasibility problem. {\sl Discrete Optimization}.  
{\bf 5} (2008), 36--52.
%Available at {\tt arXiv:math.ST/0603108}, 2008.

\bibitem{Vlach:86}
M.~Vlach. Conditions for the existence of solutions of the
three-dimensional planar transportation problem. {\sl
Disc. App. Math.} \textbf{13} (1986), 61--78.
\end{thebibliography}
\end{document}